\documentclass[12pt]{article}

\usepackage{amsmath, epsfig, cite}
\usepackage{amssymb}
\usepackage{amsfonts}
\usepackage{latexsym}
\usepackage{graphicx}
\usepackage{color}
\usepackage{ifpdf}
\usepackage{CJK}

\newtheorem{thm}{Theorem}[section]

\newtheorem{lem}[thm]{Lemma}
\newtheorem{conj}[thm]{Conjecture}

\numberwithin{equation}{section}

\makeatletter \@addtoreset{equation}{section} \makeatother

\begin{document}
\rule{0cm}{1cm}

\begin{center}
{\Large\bf Bicyclic graphs with maximal revised Szeged index }
\end{center}

\begin{center}
{\small Xueliang Li, Mengmeng Liu\\
Center for Combinatorics, LPMC\\
Nankai University, Tianjin 300071, China\\
Email:  lxl@nankai.edu.cn, liumm05@163.com}
\end{center}

\begin{center}
\begin{minipage}{120mm}
\begin{center}
{\bf Abstract}
\end{center}

{\small The revised Szeged index $Sz^*(G)$ is defined as
$Sz^*(G)=\sum_{e=uv \in E}(n_u(e)+ n_0(e)/2)(n_v(e)+ n_0(e)/2),$
where $n_u(e)$ and $n_v(e)$ are, respectively, the number of
vertices of $G$ lying closer to vertex $u$ than to vertex $v$ and
the number of vertices of $G$ lying closer to vertex $v$ than to
vertex $u$, and $n_0(e)$ is the number of vertices equidistant to
$u$ and $v$. Hansen used the AutoGraphiX and made the following
conjecture about the revised Szeged index for a connected bicyclic
graph $G$ of order $n \geq 6$:
$$
Sz^*(G)\leq \left\{
\begin{array}{ll}
(n^3+n^2-n-1)/4,& \mbox {if $n$ is odd},\\
(n^3+n^2-n)/4, & \mbox {if $n$ is even}.
\end{array}
\right.
$$
with equality if and only if $G$ is the graph obtained from the
cycle $C_{n-1}$ by duplicating a single vertex. This paper is to
give a confirmative proof to this conjecture. }

\vskip 3mm

\noindent {\bf Keywords:} Wiener index, Szeged index, Revised Szeged
index, bicyclic graph.

\vskip 3mm


\end{minipage}
\end{center}

\section{Introduction}

All graphs considered in this paper are finite, undirected and
simple. We refer the readers to \cite{bm} for terminology and
notations. Let $G$ be a connected graph with vertex set $V$ and edge
set $E$. For $u,v \in V, d(u,v)$ denotes the distance between $u$
and $v$. The {\it Wiener index} of $G$ is defined as
$$
W(G)=\displaystyle\sum_{\{u,v\}\subseteq V} d(u,v).
$$
This topological index has been extensively studied in the
mathematical literature; see, e.g., \cite{GSM,GYLL}. Let $e=uv$ be
an edge of $G$, and define three sets as follows:
$$
N_u(e) = \{w \in V: d(u,w)< d(v,w)\},
$$
$$
N_v(e) = \{w \in V: d(v,w)< d(u,w)\},
$$
$$
N_0(e) = \{w \in V: d(u,w)=d(v,w)\}.
$$
Thus, $\{N_u(e),N_v(e),N_0(e)\}$ is a partition of the vertices of
$G$ respect to $e$. The number of vertices of $N_u(e), N_v(e)$ and
$N_0(e)$ are denoted by $n_u(e), n_v(e)$ and $n_0(e)$, respectively.
A long time known property of the Wiener index is the formula
\cite{GP,W}:
$$
W(G) = \displaystyle\sum_{e=uv \in E} n_u(e) n_v(e),
$$
which is applicable for trees. Using the above formula, Gutman
\cite{G} introduced a graph invariant named the {\it Szeged index}
as an extention of the Wiener index and defined it by
$$
Sz(G) = \displaystyle\sum_{e=uv \in E}n_u(e) n_v(e).
$$
Randi\'c \cite{R} observed that the Szeged index does not take into
account the contributions of the vertices at equal distaances from
the endpoints of an edge, and so he conceived a modified version of
the Szeged index which is named the {\it revised Szeged index}. The
revised Szeged index of a connected graph $G$ is defined as
$$
Sz^*(G) = \displaystyle\sum_{e=uv \in E}\left(n_u(e)+
\frac{n_0(e)}{2}\right)\left(n_v(e)+ \frac{n_0(e)}{2}\right).
$$

Some properties and applications of this topological index have been
reported in \cite{PR,PZ}. In \cite{AH}, Aouchiche and Hansen showed
that for a connected graph $G$ of order $n$ and size $m$, an upper
bound of the revised Szeged index of $G$ is $\frac{n^2m}{4}$. In
\cite{XZ}, Xing and Zhou determined the unicyclic graphs of order
$n$ with the smallest and the largest revised Szeged indices for
$n\geq 5$, and they also determined the unicyclic graphs of order
$n$ with a unique cycle of length $r \ (3\leq r\leq n)$, with the
smallest and the largest revised Szeged indices.

In \cite{auto}, Hansen used the AutoGraphiX and made the following
conjecture:

\begin{conj} Let $G$ be a connected bicyclic graph
$G$ of order $n \geq 6$. Then
$$
Sz^*(G)\leq \left\{
\begin{array}{ll}
(n^3+n^2-n-1)/4,& \mbox {if $n$ is odd},\\
(n^3+n^2-n)/4, & \mbox {if $n$ is even}.
\end{array}
\right.
$$
with equality if and only if $G$ is the graph obtained from the
cycle $C_{n-1}$ by duplicating a single vertex (see Figure 1).
\end{conj}
It is easy to see that for bicyclic graphs, the upper bound in
Conjecture 1.1 is better than $\frac{n^2m}{4}$ for general graphs.

This paper is to give a confirmative proof to this
conjecture.

\section{Main results}

For convenience, let $B_n$ be the graph obtained from the cycle
$C_{n-1}$ by duplicating a single vertex (see Figure 1). It is easy
to check that
$$
Sz^*(B_n) =  \left\{
\begin{array}{ll}
(n^3+n^2-n-1)/4,& \mbox {if $n$ is odd},\\
(n^3+n^2-n)/4, & \mbox {if $n$ is even}.
\end{array}
\right.
$$
i.e., $B_n$ satisfies the equality of Conjecture 1.1.

So, we are left to show that for any connected bicyclic graph $G_n$
of order $n$, other than $B_n$, $Sz^*(G_n)<Sz^*(B_n)$. Using the
fact that $n_u(e)+ n_v(e)+n_0(e)=n$, we have
\begin{eqnarray*}
Sz^*(G) & = & \displaystyle \sum_{e=uv \in E} \left(n_u(e)+
\frac{n_0(e)}{2}\right)\left(n_v(e)+ \frac{n_0(e)}{2}\right) \\
&=& \displaystyle \sum_{e=uv \in E}
\left(\frac{n+n_u(e)-n_v(e)}{2}\right)
\left(\frac{n-n_u(e)+n_v(e)}{2}\right)\\
 &= & \displaystyle \sum_{e=uv \in E}
 \frac{n^2-(n_u(e)-n_v(e))^2}{4} \nonumber\\
 &= & \frac{mn^2}{4}-\frac{1}{4}\displaystyle \sum_{e=uv \in
 E}(n_u(e)-n_v(e))^2.
\end{eqnarray*}
Moreover, from $m=n+1$ we have
$$ Sz^*(G)  =
\frac{n^3+n^2}{4}-\frac{1}{4}\displaystyle \sum_{e=uv \in
 E}(n_u(e)-n_v(e))^2  \eqno(1)
$$

\setlength{\unitlength}{1mm}
\begin{center}
\begin{picture}(70,60)
\put(25,10){\circle*{1}} \put(45,10){\circle*{1}}
\put(15,20){\circle*{1}} \put(55,20){\circle*{1}}
\put(10,30){\circle*{1}} \put(20,30){\circle*{1}}
\put(15,40){\circle*{1}} \put(55,40){\circle*{1}}
\put(25,50){\circle*{1}} \put(45,50){\circle*{1}}

\put(25,10){\line(1,0){20}} \put(15,20){\line(1,-1){10}}
\put(55,20){\line(0,1){5}} \put(55,40){\line(0,-1){5}}
\put(25,50){\line(1,0){20}} \put(45,10){\line(1,1){10}}
\put(15,40){\line(1,1){10}} \put(45,50){\line(1,-1){10}}
\put(10,30){\line(1,2){5}} \put(10,30){\line(1,-2){5}}
\put(15,20){\line(1,2){5}} \put(15,40){\line(1,-2){5}}

\put(55,27.5){\circle*{0.5}} \put(55,30){\circle*{0.5}}
\put(55,32.50){\circle*{0.5}}

\put(24,3){Figure $1$: $B_n$}
\end{picture}
\end{center}

We distinguish three cases to show the conjecture. First, we
consider connected bicyclic graphs with at least one pendant edge.
Then, we consider connected bicyclic graphs without pendant edges
but with a cut vertex. Finally, we consider $2$-connected bicyclic
graphs. In the following lemmas, we deal with these cases
separately.

\begin{lem}\label{lem1}
Let $G_n$ be a connected bicyclic graph of order $n\geq 6$ with at
least one pendant edge, i.e., $\delta (G_n)=1$. Then
$$
Sz^*(G_n)<Sz^*(B_n)
$$
\end{lem}

\begin{pf}
Let $e'=xy$ be a pendant edge and $d(y)=1$. Then, for $n \geq 6,$ we
have
\begin{eqnarray*}
\displaystyle\sum_{e=uv\in E}(n_u(e)-n_v(e))^2 & \geq &
(n_x(e')-n_y(e'))^2\\
&=& (n-1-1)^2 \\
&>& n+1.
\end{eqnarray*}
Combining with equality $(1)$, the result follows.
\end{pf}
\begin{qed}
\end{qed}

\begin{lem}\label{lem2}
Let $G_n$ be a connected bicyclic graph of order $n\geq 6$ without
pendant edges but with a cut vertex, i.e., $\delta(G_n) \geq 2$ and
$\kappa(G_n)=1$. Then, we have
$$
Sz^*(G_n)<Sz^*(B_n)
$$
\end{lem}

\begin{pf}
Since $\delta(G_n)\geq 2$ and $\kappa(G_n)=1$, $G_n$ consists of two
disjoint cycles linked by a path or two cycles with a common vertex.
Assume that $C_1$ and $C_2$ are the two cycles of $G_n$, $P_t$ is
the path joining $C_1$ and $C_2$, where $t\geq 0$ is the length of
the path. Thus $|C_1|+|C_2|+t-1=n$, and $|C_1|\geq 3$ and $|C_2|\geq
3.$ Let $u \in C_1$, $v \in C_2$ be the endpoints of $P_t$. Now we
consider the four edges on the two cycles which are incident with
$u$ and $v$. Without loss of generality, we consider one of the 4
edges $e_1=uw$. Then we have
\begin{eqnarray*}
n_u(e_1)-n_w(e_1)=
n-|C_1|+\left\lfloor\frac{C_1}{2}\right\rfloor-\left\lfloor\frac{C_1}{2}\right\rfloor=n-|C_1|
\end{eqnarray*}
For the other three edges, one can get equalities similar to the
above. So we have, for $n\geq 6,$
\begin{eqnarray*}
\displaystyle\sum_{e=uv\in E}(n_u(e)-n_v(e))^2 & \geq &
2(n-|C_1|)^2+2(n-|C_2|)^2\\
&=& 2\left(2nt-2n+|C_1|^2+|C_2|^2\right)\\
&\geq & 2\left(2nt-2n+2\times \left(\frac{n+1-t}{2}\right)^2\right)\\
&=& (n-1+t)^2 \\
&>& n+1,
\end{eqnarray*}
Combining with equality $(1)$, this completes the proof.
\end{pf}
\begin{qed}
\end{qed}

For the last case, i.e., $\kappa(G_n)\geq 2$, we define a class of
graphs. A graph is called a $\Theta$-graph if it consists of three
internally disjoint paths connecting two fixed vertices. Obviously,
in this case $G_n$ must be a $\Theta$-graph.

\begin{lem}\label{lem3}
Let $G=(V,E)$ be a $\Theta$-graph, $e=uv \in E$. Then
$|n_u(e)-n_v(e)|=0$ if and only if $e$ is placed in the middle
position of an odd path of $G$.
\end{lem}

\begin{pf}
Assume that $x$ and $y$ are the vertices in $G$ with degree 3, and
$e=uv$ belongs to $P_i \ (1\leq i\leq 3)$, the $i$th path connecting
$x$ and $y$. Then, with respect to $N_u(e)$ and $N_v(e)$, there are
three cases to discuss.

\noindent {\bf Case $1$.} $x,y$ are in different sets. We claim that
$$
|n_u(e)-n_v(e)|=|b_i-a_i|,
$$
where $a_i$ (resp. $(b_i)$) is the distance between $x$ (resp. $y$)
and the edge $e$.

To see this, assume that $x \in N_u(e), \ y \in N_v(e)$. Then we
have $a_i-b_i$ vertices more in $N_u(e)$ than in $N_v(e)$ on the
path $P_i$, but on each path $P_j$ $(j\neq i)$, we have $b_i-a_i$
vertices more in $N_u(e)$ than in $N_v(e)$. Hence
$|n_u(e)-n_v(e)|=|2(b_i-a_i)+(a_i-b_i)|=|b_i-a_i|.$

\noindent{\bf Case $2$.} $x,y$ are in the same set. We claim that
$$
|n_u(e)-n_v(e)|=|V|-g,
$$
where $g$ is the length of the shortest cycle of $G$ that contains
$e$.

To see this, assume that $x, y \in N_u(e)$. Thus all vertices from
the paths $P_i$ $(j\neq i)$ are in $N_u(e)$. Therefore,
$n_v(e)=\lfloor\frac{g}{2}\rfloor$, while
$n_u(e)=\lfloor\frac{g}{2}\rfloor+|V|-g$. So
$|n_u(e)-n_v(e)|=|V|-g.$

\noindent {\bf Case $3$.} one of $x,y$ is in $N_0(e)$. We claim that
$$
|n_u(e)-n_v(e)|\geq a-1,
$$
with equality if and only if $G$ has two paths of length $a$, where
$a$ is the length of a shortest path of $G$.

To see this, assume that $x \in N_u(e)$, $y \in N_0(e)$. Then the
shortest cycle $C$ of $G$ that contains $e$ is odd. Let $z\in
V\backslash C$ be the furthest vertex from $e$ such that $z \in
N_0(e)$. Then $|n_u(e)-n_v(e)|= d(x,z)-1 \geq a+d(y,z)-1\geq a-1.$

From the above, we know that $|n_u(e)-n_v(e)|\geq 1$ in Case $2$. In
Case $3$, $|n_u(e)-n_v(e)|=0$ if $G$ has two paths of length $1$, \
which is impossible since $G$ is simple. So, $|n_u(e)-n_v(e)|=0$ if
and only if $x,y$ are in different sets and $|b_i-a_i|=0$, that is,
$e$ is placed in the middle position of an odd path of $G$.
\end{pf}
\begin{qed}
\end{qed}

Now we are ready to give our main result.
\begin{thm}
If $G_n$ is a connected bicyclic graph of order $n>6$, other than
$B_n$, then
$$
Sz^*(G_n)<Sz^*(B_n).
$$
\end{thm}
\begin{pf}
The result follows from Lemmas \ref{lem1} and \ref{lem2} for
bicyclic graphs of connectivity 1. So, we assume that $G_n$ is
$2$-connected next. Then $G_n$ must be a $\Theta$-graph. Let $x$ and
$y$ be the vertices in $G$ with degree 3, $a\leq b\leq c$ be the
lengths of the corresponding 3 paths. By Lemma \ref{lem3}, we know
that there are at most $3$ edges such that $|n_u(e)-n_v(e)|=0$. We
distinguish the following cases to proceed the proof.

\noindent {\bf Case $1$.} $3\leq a\leq b\leq c.$

Consider the six edges which are incident with $x$ and $y$. Let
$e_1=xz$ be one of them. Then, $|n_u(e)-n_v(e)|\geq 2$ from Lemma
\ref{lem3}. Similar thing is true for the other five edges. Hence
$$
\displaystyle\sum_{e=uv\in E}(n_u(e)-n_v(e))^2 \geq 2^2\times
6+(m-6-3)=m+15>m=n+1.
$$
Combining with equality $(1)$, the result follows.

\noindent {\bf Case $2$.} $2 = a < b\leq c.$

Consider the four edges which are incident with $x$ and $y$ but do
not belong to the shortest path. Let $e_1=xz$ be one of them. Then,
$|n_u(e)-n_v(e)|\geq 2$ from Lemma \ref{lem3}. Similarly, this is
true for the other three edges. Hence,
$$
\displaystyle\sum_{e=uv\in E}(n_u(e)-n_v(e))^2 \geq 2^2\times
4+(m-4-2)=m+10>m=n+1.
$$
Combining with equality $(1)$, the result follows.

\noindent {\bf Case $3$.} $1= a < b\leq c.$

If $b\geq 3$, similar to the above Case $2$, we have
$$
\displaystyle\sum_{e=uv\in E}(n_u(e)-n_v(e))^2 \geq 2^2\times
4+(m-4-3)=m+9>m=n+1.
$$
Combining with equality $(1)$, the result follows.

If $b=2$, we consider the two edges on the second longest path. Let
$e_1=xw$ be one of them. Obviously, $y \in N_0(e)$, in other words,
$|n_u(e)-n_v(e)|= d(x,z)-1 \geq a+d(y,z)-1 = d(y,z)$, where $z$ is
defined as in Case $3$ of Lemma \ref{lem3}. We claim that
$d(x,z)\geq 3$. Otherwise, if $d(x,z)\leq 2$, then $d(y,z)\leq 1$,
thus $c= d(x,z)+d(y,z)\leq 3$. It follows that $n=a+b+c-1\leq 5$, a
contradiction. Now we have
$$
\displaystyle\sum_{e=uv\in E}(n_u(e)-n_v(e))^2 \geq 2^2\times
2+(m-2-2)=m+4>m=n+1.
$$
Combining with equality $(1)$, the result follows.
\end{pf}
\begin{qed}
\end{qed}

According to our proof for Conjecture 1.1, we can also get that
among connected bicyclic graphs of order $n$, the graph
$\Theta(1,2,n-2)$ has the second-largest revised Szeged index, where
$\Theta(a,b,c)$ is a $\Theta$-graph with three paths of lengths
$a,b,c,$ respectively.

\end{document}